\documentclass[11pt]{amsart}
\usepackage{latexsym,amssymb,amsthm,amsmath,amsfonts,amsxtra,graphicx,url}
\usepackage{xypic}

\def\R{\mathbb{R}}

\def\H{\mathbb{H}}

\def\E{\mathbb{E}}
\def\S{\mathbb{S}}

\def\<{\langle}
\def\>{\rangle}

\newtheorem{theorem}{Theorem}[section]
\newtheorem{lemma}[theorem]{Lemma}
\newtheorem{proposition}[theorem]{Proposition}
\newtheorem{remark}[theorem]{Remark}

\title{On the geometry of constant angle surfaces in $Sol_3$}

\author[R. L\'opez]{Rafael L\'opez}
\address[R. L\'opez]{Departamento de Geometr\'{\i}a y Topolog\'{\i}a\\
Universidad de Granada, Spain}
\email[]{rcamino@ugr.es}

\thanks{The first author was partially supported by MEC-FEDER grant no. MTM2007-61775 and
   Junta de Andaluc\'{\i}a grant no. P06-FQM-01642}

\author[M.I. Munteanu]{Marian Ioan Munteanu}
\address[M.I. Munteanu]{University 'Al.I.Cuza' of Ia\c si\\
Faculty of Mathematics\\
Bd. Carol I, no.11\\
700506 Ia\c si\\
Romania\\
\newline
\rm
\url{http://www.math.uaic.ro/~munteanu}}
\email[]{marian.ioan.munteanu@gmail.com}

\date{\today }

\thanks{The second author was partially supported by Grant PN II ID 398/2007-2010 (Romania)}

\begin{document}



\begin{abstract}
In this paper we classify all surfaces in the 3-dimensional Lie group $Sol_3$
whose normals make constant angle with a left invariant vector field.
\end{abstract}



\maketitle

\section{Preliminaries}

The space $Sol_3$ is a simply connected homogeneous 3-dimensional manifold whose
isometry group has dimension $3$ and it is one of the eight models
of geometry of Thurston \cite{lm:Tro98}. As Riemannian manifold, the
space $Sol_3$ can be represented by $\R^3$ equipped with the metric
$$
   \widetilde g=e^{2z}dx^2+e^{-2z}dy^2+dz^2
$$
where $(x,y,z)$ are canonical coordinates of $\R^3$.
The space $Sol_3$, with the group operation
$$
  (x,y,z)\ast (x',y',z')=(x+e^{-z}x',y+e^{z}y',z+z')
$$
is a unimodular, solvable but not nilpotent Lie group and the metric $\widetilde g$ is left-invariant.
See e.g. \cite{lm:DM08,lm:Tro98}.
With respect to the metric $\tilde g$ an orthonormal basis of left-invariant vector fields is given by
$$
  e_1=e^{-z}\frac{\partial}{\partial x},\ \ e_2=e^{z}\frac{\partial}{\partial y},\ \ e_3=\frac{\partial}{\partial z}.
$$

The following transformations
$$
(x,y,z)\mapsto(y,-x,-z) \quad {\rm and}\quad (x,y,z)\mapsto(-x,y,z)
$$
span a group of isometries of $(Sol_3,g)$ having the origin as fixed point.
This group is isomorphic to the dihedral group (with $8$ elements) $D_4$.
It is, in fact, the complete group of isotropy \cite{lm:Tro98}.
The other elements of the group are
$(x,y,z)\mapsto(-x,-y,z)$, $(x,y,z)\mapsto(-y,x,-z)$, $(x,y,z)\mapsto(y,x,-z)$,
$(x,y,z)\mapsto(y,x,z)$ and $(x,y,z)\mapsto(x,-y,z)$
They can be unified as follows (cf. \cite{lm:ST05}):
$$
(x,y,x)\longmapsto (\pm e^{-c}x+a, \pm e^cy+b,z+c)
$$
$$
(x,y,z)\longmapsto(\pm e^{-c}y+a, \pm e^cx+b, z+c).
$$

It is well known that the isometry group of $Sol_3$ has dimension three.

The Levi Civita connection ${\widetilde\nabla}$ of $Sol_3$ with respect to $\{e_1,e_2,e_3\}$
is given by
$$
   \begin{array}{lll}
   {\widetilde\nabla}_{e_1} e_1=-e_3 & {\widetilde\nabla}_{e_1}e_2=0 & {\widetilde\nabla}_{e_1}e_3=\ e_1\\
   {\widetilde\nabla}_{e_2} e_1=\ 0 & {\widetilde\nabla}_{e_2}e_2=e_3 & {\widetilde\nabla}_{e_2}e_3=-e_2\\
   {\widetilde\nabla}_{e_3} e_1=\ 0 & {\widetilde\nabla}_{e_3}e_2=0 & {\widetilde\nabla}_{e_3}e_3=\ 0.
   \end{array}
$$

We recall the Gauss and Weingarten formulas

{\bf (G)} $\widetilde\nabla_XY=\nabla_XY+h(X,Y)$

{\bf (W)} $\widetilde\nabla_XN=-AX$

for every $X$ and $Y$ tangent to $M$ and for any $N$ unitary normal to $M$.
By $A$ we denote the shape operator on $M$.

\section{Constant angle surfaces in $Sol_3$ - general things}

\subsection{Motivation}

Constant angle surfaces were recently studied in product spaces ${\mathbb{Q}}_\epsilon\times\R$, where
${\mathbb{Q}}_\epsilon$ denotes the sphere $\S^2$ (when $\epsilon=+1$), the Euclidean plane $\E^2$ (when
$\epsilon=0$), respectively the hyperbolic plane $\H^2$ (when $\epsilon=-1$).
See e.g. \cite{lm:DFVV07, lm:CS07, lm:MN09, lm:DM09}.
The angle is considered between the unit normal of the surface $M$ and the tangent direction to $\R$.

It is known, for $Sol_3$, that ${\mathcal{H}}^1=\{dy\equiv0\}$ and ${\mathcal{H}}^2=\{dx\equiv0\}$
are totally geodesic foliations whose leaves are the hyperbolic plane (thought
as the upper half plane model).

On the other hand, for ${\mathbb{Q}}_\epsilon\times\R$, the foliation $\{dt\equiv0\}$ is totally geodesic too
($t$ is the global parameter on $\R$).
Trivial examples for constant angle surfaces in ${\mathbb{Q}}_\epsilon\times\R$ are furnished by totally geodesic surfaces
${\mathbb{Q}}_\epsilon\times\{t_0\}$.

Let us consider ${\mathcal{H}}^2$. It follows that the tangent plane to $\H^2$ (the leaf at each $x=x_0$)
is spanned by $\frac\partial{\partial y}$ and $\frac\partial{\partial z}$, while the unit normal is $e_1$.
So, this surface corresponds to ${\mathbb{Q}}_\epsilon\times\{t_0\}$, case in which the constant angle is 0.
Due to these reasons we give the following definition:

An oriented surface $M$, isometrically immersed in $Sol_3$, is called {\em constant
angle surface} if the angle between its normal and $e_1$ is constant in each point of the surface $M$.

\subsection{First computations}
Denote by $\theta\in[0,\pi)$ the angle between the unit normal $N$ and $e_1$.
Hence
\begin{equation*}
\widetilde g(N,e_1)=\cos\theta.
\end{equation*}

Let $T$ be the projection of $e_1$ on the tangent plane $T_pM$ of $M$ in a point $p\in M$. Thus
\begin{equation}
\label{eq:2}
e_1=T+\cos\theta N.
\end{equation}

\medskip
{\bf Case $\theta=0$.} Then $N=e_1$ and hence the surface $M$ is isometric to the hyperbolic plane
${\mathcal{H}}^2=\{dx\equiv 0\}$.

\medskip

From now on we will exclude this case.

\begin{lemma}
If $X$ is tangent to $M$ we have
\begin{itemize}
\item [1.]
      $\widetilde\nabla_Xe_1=-\widetilde g(X, e_1)e_3$,
      $\widetilde\nabla_Xe_2= \widetilde g(X,e_2)e_3$\\
      $\widetilde\nabla_Xe_3=\ \widetilde g(X, e_1)e_1 - \widetilde g(X, e_2)e_2$
\item [2.]
      $AT = -\widetilde g(N, e_3)T$,
      hence $T$ is a principal direction on the surface
\item [3.]
      $g(T,T)=\sin^2\theta$.
\end{itemize}
\end{lemma}

At this point we have to decompose also $e_2$ and $e_3$ into the tangent and the normal parts, respectively.

Let $E_1=\frac 1{\sin\theta}~T$. Consider $E_2$ tangent to $M$, orthogonal to $E_1$ and such that the basis
$\{e_1, e_2, e_3\}$ and $\{E_1, E_2,N\}$ have the same orientation. It follows that
\begin{equation}
\label{eq:7}
\left\{\begin{array}{rrcrcr}
e_1= & \sin\theta~E_1 &   & & +&\cos\theta~N\\
e_2= & \cos\alpha\cos\theta~E_1 & + & \sin\alpha~E_2 & - &\cos\alpha\sin\theta~N\\
e_3= & -\sin\alpha\cos\theta~E_1 & + & \cos\alpha~E_2 & + & \sin\alpha\sin\theta~N
\end{array}\right.
\end{equation}
and
\begin{equation}
\label{eq:16}
\left\{\begin{array}{rrcrcr}
E_1= & \sin\theta~e_1 & + & \cos\theta\cos\alpha~e_2 & - & \cos\theta\sin\alpha~e_3\\
E_2= & & & \sin\alpha~e_2 & + & \cos\alpha~e_3\\
N= & \cos\theta~e_1 & - & \sin\theta\cos\alpha~e_2 & + & \sin\theta\sin\alpha~e_3
\end{array}\right.
\end{equation}
where $\alpha$ a smooth function on $M$.

\medskip
{\bf Case $\theta=\frac\pi2$.} In this case $e_1$ is tangent to $M$ and $T=E_1$.

The metric connection on $M$ is given by
\begin{equation*}
\begin{array}{ll}
\nabla_{E_1}E_1=-\cos\alpha~E_2 & \nabla_{E_2}E_1=0\\[2mm]
\nabla_{E_1}E_2=~\cos\alpha~E_1 & \nabla_{E_2}E_2=0.
\end{array}
\end{equation*}
The second fundamental form is obtained from
\begin{equation*}
h(E_1,E_1)=-\sin\alpha~N, ~ h(E_1,E_2)=0, ~ h(E_2,E_2)=\sigma~N
\end{equation*}
where $\sigma$ is a smooth function on $M$.

Writing the Gauss formula {\bf (G)} for $X=E_1$ and $Y=E_2$, respectively for
$X=Y=E_2$ one obtains
\begin{equation*}
E_1(\alpha)=0\quad {\rm and}\quad E_2(\alpha)=\sin\alpha-\sigma.
\end{equation*}
\begin{remark}\rm
The surface $M$ is minimal if and only if $\sigma=\sin\alpha$. Since $E_1$ and $E_2$
are linearly independent, it follows that $\alpha$ is constant. Moreover, $M$ is totally
geodesic if and only if $\alpha=0$, case in which $M$ coincides with ${\mathcal{H}}^1$.
\end{remark}

Due the fact that the Lie brackets of $E_1$ and $E_2$ is $[E_1,E_2]=\cos\alpha~E_1$,
one can choose local coordinates $u$ and $v$ such that
$$
E_2=\frac\partial{\partial u}\quad{\rm and}\quad E_1=\beta(u,v)~\frac\partial{\partial v}\ .
$$
This choice implies $\alpha$ and $\beta$ fulfill the following PDE:
\begin{equation*}
\beta_u=-\beta\cos\alpha.
\end{equation*}

Since $\alpha$ depends only on $u$, it follows
$$
\beta(u,v)=\rho(v)~e^{-\int^u\cos\alpha(\tau)d\tau}
$$
where $\rho$ is a smooth function depending on $v$.

Denote by
$$\begin{array}{l}
    F:U\subset\R^2\longrightarrow M \hookrightarrow Sol_3\\[2mm]
\qquad  (u,v)\longmapsto \big(F_1(u,v),~F_2(u,v),~F_3(u,v)\big)
  \end{array}
$$
the immersion of the surface $M$ in $Sol_3$.

We have

\begin{tabbing}
(i)\qquad
$\frac\partial{\partial u}$ \= $=F_u=(F_{1,u},~F_{2,u},~F_{3,u})$\\
\> $=E_2=(\sin\alpha~e_2+\cos\alpha~e_3)_{|_{F(u,v)}}=\Big(0,~e^{F_3(u,v)}\sin\alpha,~\cos\alpha\Big)$\\[2mm]
(ii)\qquad
$\frac\partial{\partial v}$ \= $=F_v=(F_{1,v},~F_{2,v},~F_{3,v})$\\
\> $=\frac1\beta~E_1=\frac1\beta~{e_1}_{|_{F(u,v)}}=\Big(\frac1\beta~e^{-F_3(u,v)},~0,~0\Big)$.
\end{tabbing}

It follows
$$
\begin{array}{l}
F_1=F_1(v)\\
\partial_uF_2=\sin\alpha(u) e^{F_3(u,v)}\\
\partial_uF_3=\cos\alpha(u)
\end{array}
\qquad
\begin{array}{l}
\partial_vF_1=\frac1{\beta(u,v)}~e^{-F_3(u,v)}\\
F_2=F_2(u)\\
F_3=F_3(u).
\end{array}
$$

Thus we obtain
$$
\begin{array}{l}
F_1(v)=\displaystyle\int^v\frac 1{\rho(\tau)}~d\tau\\
F_2(u)=\displaystyle\int^u\Big(\sin\alpha(\tau)e^{\int^\tau\cos\alpha(s)ds}\Big)d\tau\\
F_3(u)=\displaystyle\int^u\cos\alpha(\tau)d\tau.
\end{array}
$$

Changing the $v$ parameter, one gets the following parametrization
$$
    F(u,v)=\Big(v,~\phi(u),~\chi(u)\Big)
$$
which represents a cylinder over the plane curve $\gamma(u)=\big(0,~\phi(u),~\chi(u)\big)$ where
$\phi(u)=\displaystyle\int^u\big(\sin\alpha(\tau)e^{\int^\tau\cos\alpha(s)ds}\big)d\tau$ and
$\chi(u)=\displaystyle\int^u\cos\alpha(\tau)d\tau$.
Notice that the surface is the group product between the curve
$v\mapsto(v,~0,~0)$ and the curve $\gamma$.

Let us see how the curve $\gamma$ looks like for different values of the function $\alpha$:

\begin{description}
\item [{\bf a}] $\alpha$ is a constant:
$$
\gamma(u)=\big(0,\tan\alpha~e^{u\cos\alpha},~u\cos\alpha\big)
$$

\item [{\bf b}] $\alpha(s)=s$
$$
\gamma(u)=\Big(0,\int^u\sin s~e^{\sin s} ds,~\sin u\Big)
$$

\item [{\bf c}] $\alpha(s)=s^2$
$$
\gamma(u)=\Big(0,\int^u\sin s^2~e^{\int^s\cos \tau^2d\tau}ds,~\int^u\cos s^2ds\Big)
$$

\item [{\bf d}] $\alpha(s)=\arccos(s)$, $s\in[-1,1]$
$$
\gamma(u)=\Big(0,\int^u\sqrt{1-s^2}~e^sds,~u\Big)
$$

\item [{\bf e}] $\alpha(s)=2\arctan e^{2s}$
In this case, the expression of $\gamma$ involve hypergeometric functions.
The surface $M$ is totally umbilical but not totally geodesic.

\begin{figure}[hbtp]
\begin{center}
  \includegraphics[width=55mm]{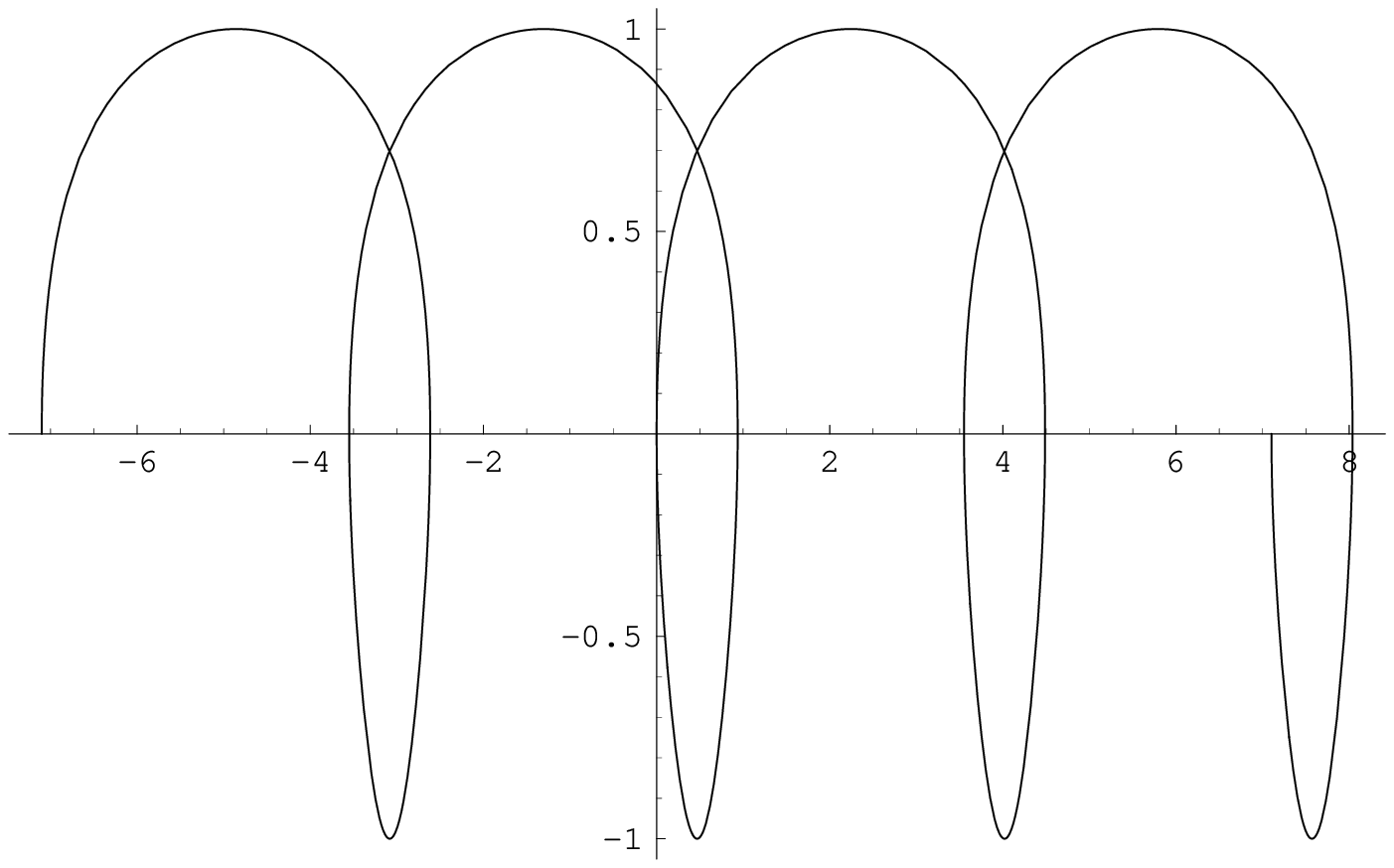} \quad
    \includegraphics[width=55mm]{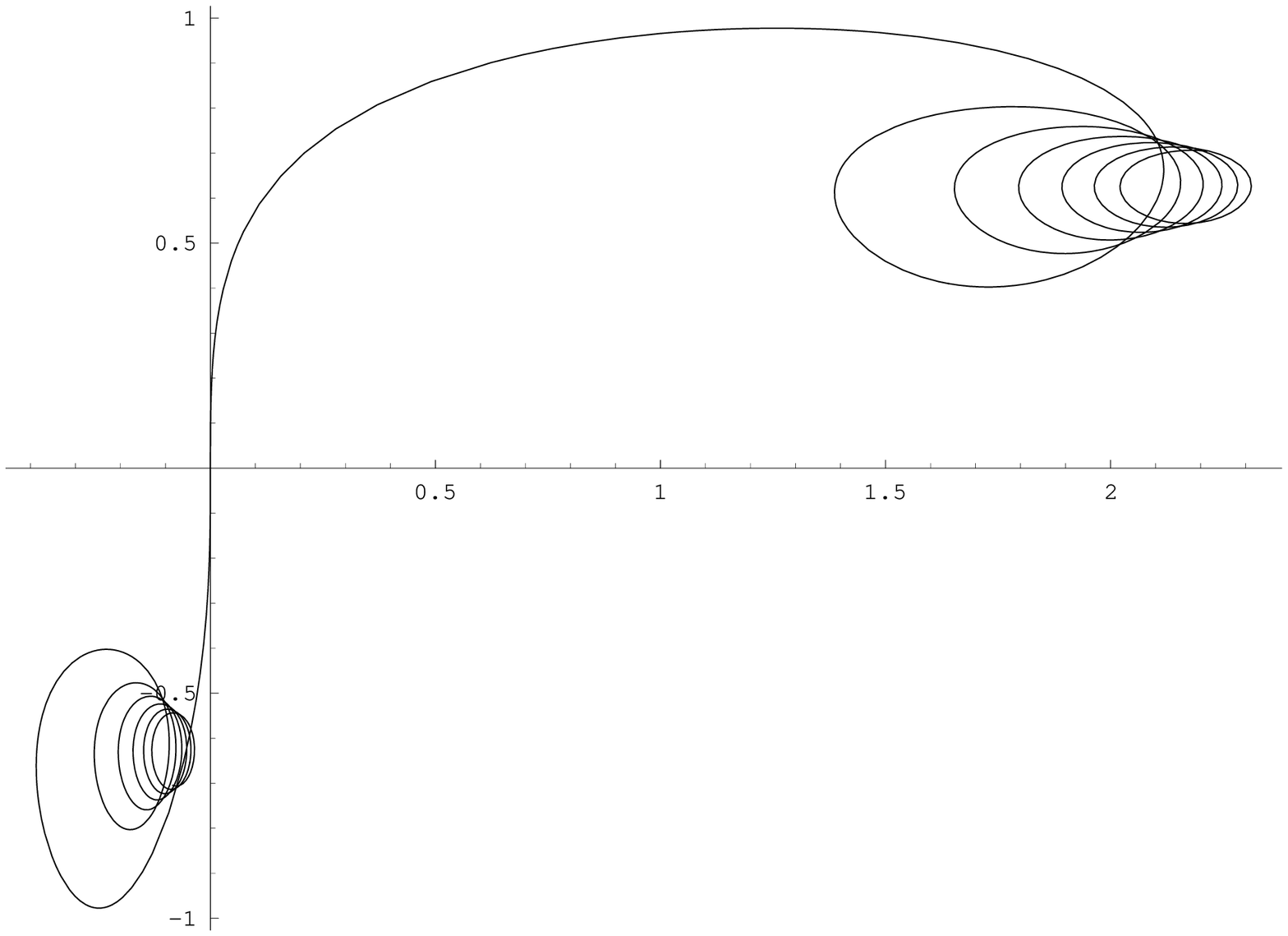}\\
\includegraphics[width=45mm]{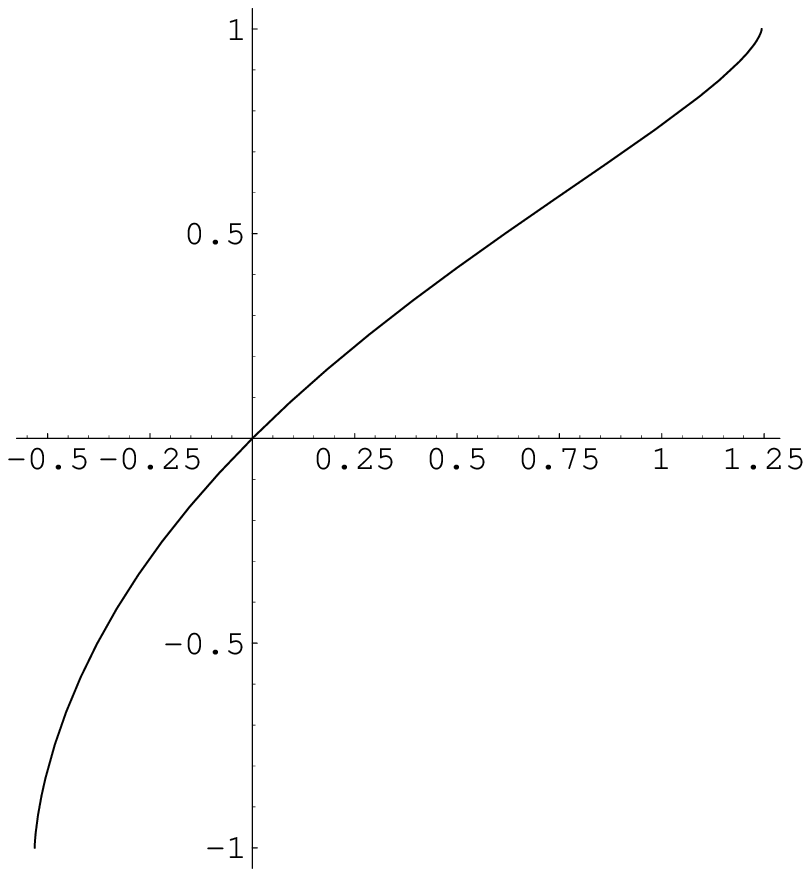} \quad
    \includegraphics[width=40mm]{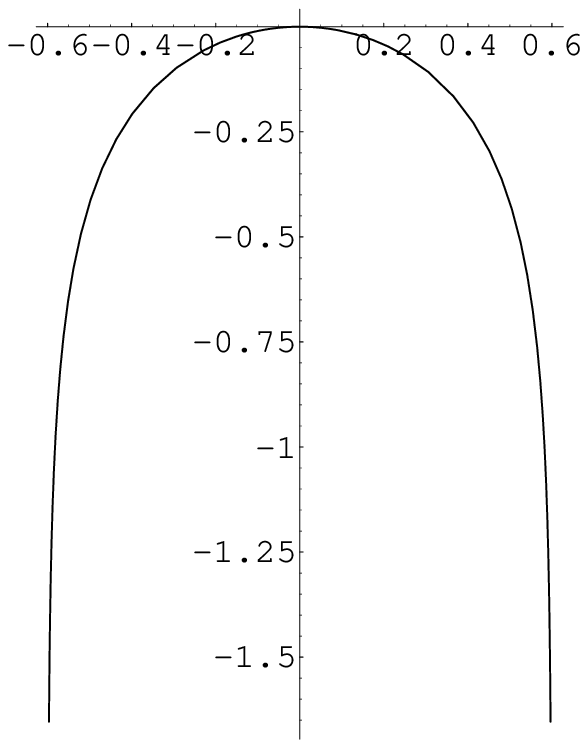}
\end{center}
\caption{Items: b, c, d and e
\label{fig:sol_03}
}
\end{figure}

\end{description}

\medskip

Coming back to the general case for $\theta$, we distinguish some particular situations for $\alpha$:

\smallskip

{\bf Case $\sin\alpha=0$.} Then $\cos\alpha=\pm1$ and the principal curvature corresponding to the
principal direction $T$ vanishes. Straightforward computations yield $\theta=\frac\pi2$ case which was
discussed before.

{\bf Case $\cos\alpha=0$.} Then $\sin\alpha=\pm1$ and the relations \eqref{eq:2} and \eqref{eq:16}
may be written in an easier way, namely, for $\sin\alpha=1$ we have
$$
\begin{array}{cl}
& e_1=\sin\theta~E_1+\cos\theta~N,\ e_2=E_2,\ e_3=-\cos\theta~E_1+\sin\theta~N\\
 & E_1=\sin\theta~e_1-\cos\theta~e_3,\ E_2=e_2,\ N=\cos\theta~e_1+\sin\theta~e_3.
\end{array}
$$
The Levi Civita connection $\nabla$ on the surface $M$ is given by
$$
\nabla_{E_1}E_1=0,\ \nabla_{E_1}E_2=0,\ \nabla_{E_2}E_1=\cos\theta~E_2,\ \nabla_{E_2}E_2=-\cos\theta~E_1.
$$
\begin{remark}
\rm Such surface is minimal.
\end{remark}
\proof
Computing the second fundamental form, one obtains
$$
h(E_1,E_1)=-\sin\theta~N,\ h(E_1,E_2)=0,\ h(E_2,E_2)=\sin\theta~N
$$
and hence the conclusion.
\endproof

In order to obtain explicit embedding equations for the surface $M$ let us choose local coordinates as follows:\\[2mm]
Let $u$ be such that $E_1=\frac\partial{\partial u}$ and $v$ such that $E_2$ and $\frac\partial{\partial v}$ are
collinear. This can be done due the fact that $[E_1,E_2]=-\cos\theta~E_2$. Considering $\frac\partial{\partial v}=b(u,v)~E_2$,
with $b$ a smooth function on $M$, since $\big[\frac\partial{\partial u}~,~\frac\partial{\partial v}\big]=0$, it follows
that $b$ satisfies $b_u-b\cos\theta=0$. This PDE has the general solution $b(u,v)=\mu(v) e^{u\cos\theta}$,
with $\mu$ a smooth function defined on certain interval in $\R$.

Denote by $F=(F_1, F_2,F_3)$ the isometric immersion of the surface $M$ in $Sol_3$.
We have
\begin{tabbing}
(i) \qquad $\frac\partial{\partial u}$ \= $=F_u=\big(\partial_uF_1,~\partial_uF_2,~\partial_uF_3\big)$\\[2mm]
    \> $=E_1=\sin\theta~{e_1}_{|_{F(u,v)}}-\cos\theta~{e_3}_{|_{F(u,v)}}$
       $=\big(\sin\theta e^{-F_3(u,v)},~0,~-\cos\theta\big)$\\[3mm]
(ii) \qquad $\frac\partial{\partial v}$ \= $=F_v=\big(\partial_vF_1,~\partial_vF_2,~\partial_vF_3\big)$\\[2mm]
    \> $=\mu(v) e^{u\cos\theta}~E_2$ $=\mu(v)e^{u\cos\theta}~{e_2}_{|_{F(u,v)}}$\\[2mm]
    \> $=\big(0,~\mu(v) e^{u\cos\theta+F_3(u,v)},~0\big)$.
\end{tabbing}

Looking at (i) we immediately get
\begin{itemize}
\item [$\bullet$] the third component:
$
   F_3(u,v)=-u\cos\theta+\zeta(v)\ , \quad \zeta\in C^\infty(M)
$
\item [$\bullet$] the second component:
$ F_2(u,v)=F_2(v) $.
\end{itemize}

Replacing in (ii) we obtain
\begin{itemize}
\item[$\bullet$] the third component: $\zeta(v)=\zeta_0\in\R$
\item[$\bullet$] the second component: $F_2(v)= e^{~\zeta_0}\displaystyle\int^v \mu(\tau) {\rm d}\tau$
\item[$\bullet$] the first component: $F_1(u,v)=F_1(u)$.
\end{itemize}

Going back in (i) and taking the first component one gets
$$
F_1(u)=e^{-\zeta_0} \tan\theta e^{u\cos\theta}+{\rm constant}.
$$
Since the map $(x,y,z)\longmapsto(x+c,y,z)$ is an isometry for $Sol_3$, we can take the previous constant to be $0$.
Moreover, the map $(x,y,z)\longmapsto(e^{-c}x,e^cy,z+c)$ is also an isometry of the ambient space, so $\zeta_0$ may be assumed
to be also $0$.

Consequently, one obtains the following parametrization for the surface $M$
$$
F(u,v)=\big(\tan\theta~e^{u\cos\theta},~\int^v\mu(\tau) {\rm d}\tau,~-u\cos\theta\big).
$$

Finally, we can change the parameter $v$ such that $\mu(v)=1$. One can state the following
\begin{proposition}
\label{prop2.4}
The surface $M$ given by the parametrization
\begin{equation}
\label{eq:F11}
F(u,v)=\left(\tan\theta~e^{u\cos\theta},~v,~-u\cos\theta\right)
\end{equation}
is a constant angle surface in $Sol_3$.
\end{proposition}
Notice that this surface is a (group) product between the curve $v\mapsto(0,~v,~0)$ and the plane curve
$\gamma(u)=(\tan\theta~e^{u\cos\theta},~0,~-u\cos\theta)$.

The angle $\theta$ is an arbitrary constant.
Moreover, the curvature of $M$ is a negative constant $-\cos^2\theta$.
Analogue results are obtained if $\cos\alpha=-1$.

\bigskip

From now on we will deal with $\alpha$ and $\theta$ different from the situations above.

\begin{lemma}The Levi Civita connection $\nabla$ on $M$ and the second fundamental form $h$ are given by
\begin{equation}
\label{eq:LC}
\left\{
\begin{array}{ll}
\nabla_{E_1}E_1=-\cos\alpha~E_2, &
\nabla_{E_1}E_2=\cos\alpha~E_1\\[2mm]
\nabla_{E_2}E_1=\sigma\cot\theta~E_2, &
\nabla_{E_2}E_2=-\sigma\cot\theta~E_1
\end{array}
\right.
\end{equation}
\begin{equation}
\label{eq:second_ff}
h(E_1,E_1)=-\sin\theta\sin\alpha~N,\quad
h(E_1,E_2)=0,\quad
h(E_2,E_2)=\sigma~N.
\end{equation}
The matrix of the Weingarten operator $A$ with respect to the basis $\{E_1, E_2\}$ has the
following expression
$$
        A=\left(\begin{array}{cc}-\sin\alpha\sin\theta & 0\\ 0 & \sigma
          \end{array}\right)
$$
for a certain function $\sigma\in C^\infty(M)$.
\end{lemma}
Moreover, the Gauss formula yields
\begin{subequations}
\renewcommand{\theequation}{\theparentequation .\alph{equation}}
\label{eq:alpha_cond}
\begin{eqnarray}
\label{mn:alpha_cond_E1}
&& E_1(\alpha)=2\cos\theta\cos\alpha\\
\label{mn:alpha_cond_E2}
&& E_2(\alpha)=\sin\alpha-\frac\sigma{\sin\theta}
\end{eqnarray}
\end{subequations}
and the compatibility condition
$$
\left(\nabla_{E_1}E_2-\nabla_{E_2}E_1\right)(\alpha)=[E_1,E_2](\alpha)=E_1(E_2(\alpha))-E_2(E_1(\alpha))
$$
gives rise to the following differential equation
\begin{equation}
\label{eq:PDE_sigma}
E_1(\sigma)+\sigma\cos\theta\sin\alpha+\sigma^2\cot\theta=2\sin\theta\cos\theta\sin^2\alpha.
\end{equation}

\begin{remark}
\rm The curvature of $M$ is equal to $2\sin^2\alpha\sin^2\theta-\sigma\sin\alpha\sin\theta-1$.
\end{remark}

We are looking for a coordinate system $(u,v)$ in order to determine the embedding equations of the surface.
Let us take the coordinate $u$ such that $\frac\partial{\partial u}=E_1$.
Concerning $v$, we will discuss later about it.

Let point our attention on \eqref{mn:alpha_cond_E1} which can be re-written as
$$
\partial_u\alpha=2\cos\theta\cos\alpha.
$$
Solving this PDE one gets
$$
\sin\alpha=\tanh(2u\cos\theta+\psi(v))
$$
where $\psi$ is a smooth function on $M$ depending on $v$.
Notice that, apparently the equation has also a second solution
$\sin\alpha=\coth(2u\cos\theta+\psi(v))$. This is not valid because
$\coth$ takes values in $(-\infty,-1)$ or in $(1,+\infty)$.

Now, let us take $v$ in such way that $\frac{\partial\alpha}{\partial v}=0$, namely $\psi$ is a constant,
denote it by $\psi_0$. It follows that $\alpha$ is given by
\begin{equation}
\label{eq:sol_alpha}
\sin\alpha=\tanh(\bar u)
\end{equation}
where $\bar u=2u\cos\theta+\psi_0$.

At this point, the equation \eqref{eq:PDE_sigma} becomes
\begin{equation}
\label{eq:sigma_u}
\sigma_u+\cot\theta~\big(\sigma+2\sin\alpha\sin\theta\big)\big(\sigma-\sin\alpha\sin\theta\big)=0.
\end{equation}

Since $\frac\partial{\partial v}$ is tangent to $M$, it can be decomposed in the basis $\{E_1,E_2\}$.
Thus, there exist functions $a=a(u,v)$ and $b=b(u,v)$ such that
$$
\frac\partial{\partial v}=aE_1+bE_2.
$$

Due to the choice of the coordinate $v$ we have
$$
0=\frac{\partial \alpha}{\partial v}=a\cdot2\cos\theta\cos\alpha+
  b\left(\sin\alpha-\frac\sigma{\sin\theta}\right).
$$

{\bf a.} The case $b=0$ implies $\cos\theta=0$ or $\cos\alpha=0$.
Both situations were studied separately.

{\bf b.} Consider $b\neq0$. Let us denote by $p(u,v)=\frac ab\ $.
Hence the equality above yields
\begin{equation}
\label{eq:30}
\sigma=\sin\theta\sin\alpha + p~\sin2\theta\cos\alpha.
\end{equation}

On the other hand
$$
0=\left[\frac\partial{\partial u}~,~\frac\partial{\partial v}\right]=a_uE_1+b_uE_2+
   b\big(\cos\alpha E_1-\sigma\cot\theta E_2\big).
$$
Hence
\begin{equation}
\label{eq:31}
\left\{
\begin{array}{l}
a_u+b\cos\alpha=0\\[2mm]
b_u-b\sigma\cot\theta=0.
\end{array}\right.
\end{equation}

If we take in \eqref{eq:30} the derivative with respect to $u$, and combining with \eqref{eq:sigma_u}, it follows
\begin{equation}
\label{eq:PDE_p}
p_u+\cos\alpha+p\cos\theta\sin\alpha+2p^2\cos^2\theta\cos\alpha=0.
\end{equation}
Straightforward computations yield the general solution for this
equation (see the Appendix), namely
\begin{equation}
\label{eq:sol_p}
    p(u,v)=\pm \frac1{\cos\theta\sinh\bar u+\varepsilon~\frac{\cosh^{\frac32}\bar u}{-I(u)+\Lambda(v)}}
\end{equation}
where $\varepsilon=0,1$ and $\Lambda$ is a certain function depending on $v$.


Let $F:U\subset\R^2\longrightarrow M\hookrightarrow Sol_3$,
$(u,v)\longmapsto\big(F_1(u,v),~F_2(u,v),~F_3(u,v)\big)$
be the immersion of the surface $M$ in $Sol_3$.
We have

\begin{tabbing}
{\bf I.}
$\partial_u$ \= $=F_u=(F_{1,u},~F_{2,u},~F_{3,u})$\\
\> $=E_1$ \= $=\sin\theta~{e_1}_{|_{F(u,v)}}+\cos\theta~{e_2}_{|_{F(u,v)}}-
   \cos\theta\sin\alpha~{e_3}_{|_{F(u,v)}}$\\
\> \> $=\left(\sin\theta~e^{-F_3(u,v)},~\cos\theta\cos\alpha~e^{F_3(u,v)},-\cos\theta\sin\alpha\right)$
\end{tabbing}
which implies
\begin{subequations}
\renewcommand{\theequation}{\theparentequation .\alph{equation}}
\begin{eqnarray}
\label{eq:M1}
&&\partial_uF_1=\sin\theta~e^{-F_3(u,v)}
\end{eqnarray}
\begin{eqnarray}
\label{eq:M2}
&&\partial_uF_2=\cos\theta\cos\alpha~e^{F_3(u,v)}
\end{eqnarray}
\begin{eqnarray}
\label{eq:M3}
&&\partial_uF_3=-\cos\theta\sin\alpha.
\end{eqnarray}
\end{subequations}
From the last equation one immediately obtains
\begin{equation}
\label{eq:M4}
F_3(u,v)=-\frac12\log\cosh(\bar u)+\zeta(v)
\end{equation}
where $\zeta$ is a smooth function. Replacing this expression in
\eqref{eq:M1} and \eqref{eq:M2}, one gets
\begin{equation}
\label{eq:M5}
F_1=\sin\theta~e^{-\zeta(v)}(I(u)+f_1(v))
\end{equation}
\begin{equation}
\label{eq:M6}
F_2=\pm\cos\theta~e^{\zeta(v)}(J(u)+f_2(v))
\end{equation}
where $\displaystyle I(u)=\int\limits^u\sqrt{\cosh(2\tau\cos\theta+\psi_0)}d\tau$,
$\displaystyle J(u)=\int\limits^u\cosh^{-\frac32}(2\tau\cos\theta+\psi_0)d\tau$ and
$f_1$, $f_2$ are some smooth functions which will be determined
in what follows.

\bigskip

\begin{tabbing}
{\bf II.}
$\partial_v$ \= $=F_v$ \= $=(F_{1,v},~F_{2,v},~F_{3,v})$\\
\> $=a(u,v)E_1+b(u,v)E_2$\\
\> $=a(u,v)\big(\sin\theta~{e_1}_{|_{F(u,v)}}+\cos\theta\cos\alpha~{e_2}_{|_{F(u,v)}}-
    \cos\theta\sin\alpha~{e_3}_{|_{F(u,v)}}\big)+$\\
\>\> $+b(u,v)\big(\sin\alpha~{e_2}_{|_{F(u,v)}}+\cos\alpha~{e_3}_{|_{F(u,v)}}\big)$.
\end{tabbing}

It follows
\begin{subequations}
\renewcommand{\theequation}{\theparentequation .\alph{equation}}
\begin{eqnarray}
\label{eq:M7}
&&\partial_vF_1=a(u,v)\sin\theta~e^{-F_3(u,v)}
\end{eqnarray}
\begin{eqnarray}
\label{eq:M8}
&&\partial_vF_2=\big(a(u,v)\cos\theta\cos\alpha+b(u,v)\sin\alpha\big)~e^{F_3(u,v)}
\end{eqnarray}
\begin{eqnarray}
\label{eq:M9}
&&\partial_vF_3=-a(u,v)\cos\theta\sin\alpha+b(u,v)\cos\alpha.
\end{eqnarray}
\end{subequations}

From \eqref{eq:M4} and \eqref{eq:M9} we have
$$
     -a(u,v)\cos\theta\sin\alpha+b(u,v)\cos\alpha=\zeta'(v)
$$
and from \eqref{eq:M5} and \eqref{eq:M7} we obtain
\begin{equation}
\label{eq:M11}
\zeta'(v)\big(I(u)+f_1(v)\big)-f_1'(v)+a(u,v)\sqrt{\cosh(\bar u)}=0.
\end{equation}
Taking the derivative with respect to $u$, one gets
\begin{equation}
\label{eq:M12}
\zeta'(v)+a_u(u,v)+a(u,v)\cos\theta\tanh(\bar u)=0.
\end{equation}

The equation in $a$ has the solution
\begin{equation}
\label{eq:M13}
a(u,v)=\frac{-\zeta'(v)I(u)+\xi(v)}{\sqrt{\cosh(\bar u)}}
\end{equation}
\begin{equation}
\label{eq:M15}
b(u,v)=\pm\left[\frac{\cos\theta\sinh(\bar u)}{\sqrt{\cosh(\bar u)}}
     \big(-\zeta'(v)I(u)+\xi(v)\big)+\zeta'(v)\cosh(\bar u)\right].
\end{equation}

\medskip

Recall that $p(u,v)=\frac{a(u,v)}{b(u,v)}$. We immediately notice that
the general solution given by \eqref{eq:sol_p} is obtained with the following identification:
$\varepsilon=0\iff\zeta'(v)=0$ and $\varepsilon=1\iff\Lambda(v)=\frac{\xi(v)}{\zeta'(v)}$.
It follows
$$
   p(u,v)=\pm\frac{1}{\cos\theta\sinh(\bar u)+\frac{\zeta'(v)\cosh^{\frac32}\bar u}{-\zeta'(v)I(u)+\xi(v)}}~.
$$

At this point we will obtain the parametrization of the surface in the following way.

{\bf 1.} Combining \eqref{eq:M13} with \eqref{eq:M11} one gets
$f_1'(v)-\zeta'(v)f_1(v)-\xi(v)=0$ which has the solution
$f_1(v)=e^{\zeta(v)}\displaystyle\int\limits^v\xi(\tau)e^{-\zeta(\tau)}d\tau$. Thus
$$
    F_1(u,v)=\sin\theta\Big(e^{-\zeta(v)}I(u)+\int\limits^v\xi(\tau)e^{-\zeta(\tau)}d\tau\Big).
$$

{\bf 2.} Similarly, replace \eqref{eq:M6} in \eqref{eq:M8} one obtains
\begin{equation}
\label{eq:M18}
\begin{array}{l}
\cos\theta \big(f_2'(v)+\zeta'(v)f_2(v)\big)+\qquad\\[2mm]
\qquad
   +\zeta'(v)\big(\cos\theta (I(u)+J(u))-
       \frac{\sinh(\bar u)}{\sqrt{\cosh(\bar u)}}\big)=\cos\theta~\xi(v).
\end{array}
\end{equation}

We have
$$
   a(u,v)\cos\theta\cos\alpha+b(u,v)\sin\alpha=\pm\zeta'(v)
      \big(\sinh(\bar u)-\cos\theta~I(u)\sqrt{\cosh(\bar u)}\big)
$$
and
$$
   \cos\theta\big(I(u)+J(u)\big)-
   \frac{\sinh(\bar u)}{\sqrt{\cosh(\bar u)}}={\rm constant}
$$
which can be incorporated in the primitives $I(u)$ or $J(u)$.
It follows that $f_2$ satisfies the following ODE
$f_2'(v)+\zeta'(v)f_2(v)=\xi(v)$
which has the solution $f_2(v)=e^{-\zeta(v)}\displaystyle\int\limits^v\xi(\tau)e^{\zeta(\tau)}d\tau$.
Thus
$$
F_2(u,v)=\pm\cos\theta\Big(e^{\zeta(v)} J(u)+\int\limits^v\xi(\tau)e^{\zeta(\tau)}d\tau\Big).
$$

We conclude with the following result
\begin{theorem}
A general constant angle surface in $Sol_3$ can be parameterized as
\begin{equation}
F(u,v)=\gamma_1(v)*\gamma_2(u)
\end{equation}
where
\begin{subequations}
\renewcommand{\theequation}{\theparentequation .\alph{equation}}
\begin{eqnarray}
\label{eq:gamma1}
&&\gamma_1(v)=\Big(\sin\theta\int\limits^v\xi(\tau)e^{-\zeta(\tau)}d\tau,~
    \pm\cos\theta\int\limits^v\xi(\tau)e^{\zeta(\tau)}d\tau,~\zeta(v)\Big)
\end{eqnarray}
\begin{eqnarray}
\label{eq:gamma2}
&&\gamma_2(u)=\Big(\sin\theta ~I(u),~\pm\cos\theta ~J(u),~-\frac12\log\cosh\bar u\Big)
\end{eqnarray}
\end{subequations}
and $\zeta$, $\xi$ are arbitrary functions depending on $v$.
\end{theorem}
The curve $\gamma_2$ is parametrized by arclength.

\begin{remark}\rm
The only minimal constant angle surfaces in $Sol_3$ are: {\bf (i)} the hyperbolic plane ${\mathcal{H}}^2$
(for $\theta=0$); {\bf (ii)} the hyperbolic plane ${\mathcal{H}}^1$ (for $\theta=\frac\pi2$); {\bf (iii)}
surfaces furnished by Proposition~\ref{prop2.4}.
\end{remark}
\proof
In the general case when $\theta$ is different from $0$ and $\frac\pi2$ and $\alpha$ is
such that $\sin\alpha$ and $\cos\alpha$ do not vanish, the minimality condition can be written as
$\sigma=\sin\alpha\sin\theta$. But this relation is impossible due to \eqref{eq:30} and \eqref{eq:PDE_p}.
\endproof

\medskip

{\bf Final Remark.}
In order to define constant angle surfaces in $Sol_3$ we have considered $e_1$ as the
direction with which the normal to the surface makes constant angle. Since both
$\mathcal{H}^1$ and $\mathcal{H}^2$ are totally geodesic foliations one can also propose
$e_2$ as a candidate to the preferred direction. If this is the choice, one can define
constant angle surfaces in $Sol_3$ to be those surfaces $M$ whose unit normals make constant
angle with $e_2$ in each point of $M$. Analogue computations give rise to similar results.
Since the differences are insignificant we do not give any detail for this problem.

\section{Appendix: Solution of PDE}

{\bf Problem.} Solve the equation
$p_u+\cos\alpha+p\cos\theta\sin\alpha+2p^2\cos^2\theta\cos\alpha=0 $.

{\bf Solution.} Denote by $\bar u=2u\cos\theta+\psi_0$.

Let $q:=\frac 1p$; it follows that $q$ satisfies
$$
q_u-q^2\cos\alpha-q\cos\theta\sin\alpha-2\cos^2\theta\cos\alpha=0.
$$

Let $A:=q-\cos\theta\sinh\bar u$. It follows
$q_u=A_u+2\cos^2\theta\cosh\bar u$. Hence, $A$ satisfies
$$
A_u-3A\cos\theta\sinh\bar u-\frac 1{\cosh\bar u}~A^2=0.
$$

Let $B:=A\cosh^{-\frac32}\bar u$. It follows
$A_u=3B\cos\theta\sinh\bar u\cosh^{\frac12}\bar u+B_u\cosh^{\frac32}\bar u$.
Thus, $B$ satisfies
$$
B_u-B^2\cosh^{\frac12}\bar u=0.
$$

Hence either $B=0$ or $\frac1{B(u,v)}=-I(u)+\Lambda(v)$, for a smooth $\Lambda$.

\medskip

If $B=0$ then $A=0$, $q=\cos\theta\sinh\bar u$.

\qquad
$q\neq0$ if and only if $\theta\neq\frac\pi2$ and $\bar u\neq0$.

One gets
$$
p=\frac1{\cos\theta\sinh\bar u}~.
$$

If $B\neq0$ then
$$
    q(u,v)=\cos\theta\sinh\bar u+\frac{\cosh^{\frac32}\bar u}{-I(u)+\Lambda(v)}.
$$
These solutions correspond to
{\bf 1.} $\zeta'=0$ and {\bf 2.} $\Lambda(v)=\frac{\xi(v)}{\zeta'(v)}$

\bigskip


\end{document}